\documentstyle{amsart}

\def\ra{\rightarrow}

\def\Z{\mbox{$\Bbb Z$}}

\def\ph{\varphi}

\def\and{\hbox{ and }}
\def\s#1{\mbox{$(-1)^{#1}$}}

\def\ip#1#2{\left<#1,#2\right>}

\def\brf{\left\{\cdot,\cdot\right\}}

\def\blacksquare{\quad \vrule height 6pt width 6pt}

\newtheorem{thm}{Theorem}

\def\del{\partial}

\def\ie{\hbox{\it i.e.}}
\def\ainf{\hbox{$A_\infty$}}
\def\linf{\mbox{$L_\infty$}}
\def\zt{\hbox{$\Z_2$}}
\def\ztz{\hbox{$\Z_2\times\Z$}}

\def\tph{\hbox{$\tilde \ph$}}

\def\td{\hbox{$\tilde d$}}

\def\hom{\mbox{\rm Hom}}

\def\tns{\otimes}
\def\mtns{\tns\cdots\tns}
\def\modot{\odot\cdots\odot}

\def\mwedge{\wedge\cdots\wedge}
\def\modot{\odot\cdots\odot}
\def\mplus{+\cdots+}
\def\mcom{,\cdots,}
\def\iso{\kern.35em{\raise3pt\hbox{$\sim$}\kern-1.1em\to}
         \kern.3em}

\def\sg{\hbox{$\sigma_\Gamma$}}

\def\go{\hbox{$\Gamma_{\rm or}$}}
\def\gop{\hbox{${\Gamma^{'}}_{\rm or}$}}

\def\k{\mbox{\bf k}}

\def\shf{\operatorname{Sh}}
\def\sh#1#2{\shf(#1,#2)}

\def\s#1{\mbox{$(-1)^{#1}$}}

\def\ip#1#2{\left<#1,#2\right>}
\def\blacksquare{\quad \vrule height 6pt width 6pt}

\begin{document}
\nocite{sta3,lod,conn,hoch}
\author{Michael Penkava}
\thanks{Partially supported by NSF grant DMS-94-0411. 
The author would also like
to thank the University of Washington, Seattle for hosting him this 
summer.}
\address{University of California\\
Davis, CA 95616}
\email{michae@@math.ucdavis.edu}
\title{Infinity Algebras and the Homology of Graph Complexes}
\maketitle
\input{epsf}
\begin{abstract}
An \linf\ algebra is a generalization of a Lie algebra
\cite{ls,lm,ss}. Given an
\linf\ algebra with an invariant inner product, we construct a cycle in
the homology of the complex of metric ordinary graphs. Since the cyclic
cohomology of a Lie algebra determines infinitesimal deformations
of the algebra into an \linf\ algebra, this construction shows that 
a cyclic cocycle of a Lie algebra determines a cycle in the 
homology of the graph complex. This result was suggested by a remark
by M. Kontsevich in \cite{kon3} that every finite dimensional Lie
algebra with an invariant inner product determines a cycle in the
graph complex. In a joint paper with A. Schwarz 
\cite{ps2}, we proved that an
\ainf\ algebra with an invariant inner product determines a cycle in
the homology of the complex of metric ribbon graphs. In this article,
a simpler proof of this fact is given. Both constructions are based
on the ideas presented in \cite{pen1}.
\end{abstract}
\section{Introduction}
In \cite{pen1}, a definition of the cohomology of an \linf\
algebra was given, and it was shown that the cohomology of a Lie 
algebra classifies the (infinitesimal)
deformations of the Lie algebra into an \linf\
algebra. The cyclic cohomology of a Lie algebra was defined, and it
was shown that the cyclic cohomology of a Lie algebra with an invariant
inner product classifies the (infinitesimal) deformations of the
Lie algebra into an \linf\ algebra preserving the invariant inner product.

In a previous joint article with Albert Schwarz \cite {ps2}, similar results
were obtained about the cohomology and cyclic cohomology of an associative
algebra classifying the deformations of the algebra into an \ainf\ algebra.
In that article, we defined Hochschild cohomology and cyclic cohomology
of an \ainf\ algebra, and also applied these considerations to construct
cycles in the homology of the complex of metric ribbon graphs, which is
closely associated to the homology of the moduli spaces of Riemann surfaces
with fixed genus and number of marked points. 

The cohomology of \ainf\
algebras was studied in more detail in \cite{pen1}, from the point of
view of the coalgebra structure on the tensor coalgebra of a vector space,
and the coalgebra structure on its parity reversion (we assume the space
in question has a \zt-grading). A bracket structure, the
Gerstenhaber bracket, was defined on the
tensor algebra. 
This bracket arises from the bracket of coderivations on the tensor
coalgebra of the parity reversion. The differential in the cohomology
of an \ainf\ algebra is given by the Gerstenhaber bracket with the
cochain determined by the \ainf\ structure, much as the differential of 
cochains of an associative algebra is given by the Gerstenhaber bracket
with the cochain given by the associative multiplication, as was shown
in \cite{gers}.

When an \ainf\ algebra possesses a (graded symmetric) inner product, there is
a notion of a cyclic cochain. The property of cyclicity is preserved 
by the Gerstenhaber bracket, and the differential is given by the
bracket when the inner product is invariant. One can also define 
a notion of cyclicity in the tensor coalgebra of the parity reversion,
by means of the isomorphism that exists between the two structures.
We will give this definition here. 
Another important fact is that
symmetric inner products on a vector space give rise to graded
antisymmetric inner products on its parity reversion. 

The cohomology of an \linf\ algebra is given by considering the
exterior coalgebra of a vector space. The parity reversion of the 
exterior coalgebra is the symmetric coalgebra of the parity reversion.
There is a bracket defined on the exterior coalgebra which arises
from the bracket of coderivations of the symmetric coalgebra of the
parity reversion, and the differential in the cohomology of an \linf\
algebra is given by the bracket in much the same manner as in the \ainf\
case.

The notion of a cyclic cochain for an \linf\ algebra with an invariant
graded symmetric inner product yields a theory of cyclic cohomology
of \linf\ algebras, where again, the differential is given by a bracket.
We will need to consider here how the 
defining relations for an \linf\ algebra are expressed in
the corresponding theory on the
symmetric coalgebra of the parity reversion, which will possess a
graded antisymmetric inner product.

In \cite{kon3}, an orientation on the complex of metric ordinary graphs
was given. We will need an equivalent formulation of this orientation.
Similarly, in \cite{ps2} we gave a definition (following 
Kontsevich) of the orientation on the complex of metric ribbon graphs;
we shall also need to reformulate this definition.

We will explain the graph complexes and the definition of homology in
some detail, and then formulate the definitions of the cycles we are
interested in. The results in the next two sections are really just
a reformulation of results which were proved in \cite{ps2,pen1}.
\section{\ainf\ algebras}\label{ainfsect}
In what follows we shall assume that $V$ is a 
finite dimensional \zt-graded vector space over
a field \k\ of characteristic zero. 
The \zt-grading means that $V$ has a preferred decomposition 
$V=V_0\bigoplus V_1$. An element in $V_i$ has even parity, when $i=0$,
and odd parity when $i=1$.
The parity reversion $\Pi V$ is the same space, but with reversed parity;
\ie, $(\Pi V)_0=V_1$ and $(\Pi V)_1=V)0$. 
The tensor coalgebra $T(V)$ is defined by $T(V)=\bigoplus_{n=1}^\infty V^n$,
where $V^n$ is the $n$-th tensor power of $V$. An \ainf\ algebra
structure on $V$ is given by an element (cochain)
 $m\in C(V)=\hom(T(V),V)$, in other
words, a collection of maps $m_k:V^k\ra V$, satisfying,  for any
$v_1\mcom v_n\in V$, the relations
\begin{equation}\label{mysigns}
\sum
\begin{Sb}
k+l=n+1
\\
0\le i\le n-k\\
\end{Sb}
\s{s_{i,k,n}}
m_l(v_1\mcom v_i,m_k(v_{i+1}\mcom v_{i+k}),v_{i+k+1}\mcom v_{n})=0.
\end{equation}
where $s_{i,k,n}=(v_1\mplus v_i)k+i(k-1)+n-k$.
In addition, we require that $m_k$ have parity $k$, so that $m_2$ is an
even map, and $m_1$ is an odd map. The map $m_1$ determines a differential
on $V$, and is a derivation with respect to the product that $m_2$ 
determines on $V$. Also, the map $m_2$ is associative up to a homotopy
determined by the element $m_3$. These statements follow from the
relations given in equation (\ref{mysigns}) for $n=1,$ 2, and 3. The
last fact  accounts for the term {\bf strongly homotopy associative 
algebras}, which is another name for \ainf\ algebras.

If we let $\Pi W=V$, then there is a natural isomorphism $\eta$ between
$T(W)$ and $T(V)$, given by
\begin{equation}
\eta(w_1\mtns w_n)=\s{(n-1)v_1+(n-2)v_2\mplus (1)v_{n-1}}\pi w_1\mtns
\pi w_n.
\end{equation}
Note that the  restriction of $\eta$ to $W^k$ is a map
$\eta_k:W^k\ra V^k$, which has parity $k$, so this map is unusual
because it is not homogeneous.
The map $\eta$ induces an isomorphism between $C(V)$ and $C(W)$, given by 
$m_k\mapsto d_k$, with $d_k=\eta^{-1}\circ m_k\circ \eta$. 
Now any element(cochain) $d\in C(W)$ extends to a 
\zt-graded coderivation on the coalgebra
$T(W)$.  In \cite{pen1}, it was shown that $m$ determines an \ainf\ 
structure on $V$ precisely when $d=\eta^{-1}\circ m\circ \eta$ 
determines an odd codifferential on $T(W)$. (Actually this fact is well
known, and really is the basis for the definition of \ainf\ algebra
in the first place.) We shall regard $W$, equipped with the structure
given by $d$, as an \ainf\ algebra.
Now the condition that $d$ determines a codifferential 
on $C(W)$ is given by
\begin{equation}
\sum
\begin{Sb}
k+l=n+1\\
\\
0\le i\le n-k
\end{Sb}
\s{(w_1\mplus w_i)d_k}
d_l(w_1\mcom w_i,d_k(w_{i+1}\mcom w_{i+k}),w_{i+k+1}\mcom w_n)=0,
\end{equation}
for any $w_1\mcom w_n\in W$.
Because of the realization of elements in $C(W)$ as coderivations, there
is a natural \zt-graded bracket on $C(W)$. The bracket originates with
the associative product on $C(W)$ given by
\begin{multline}
d\delta(w_1\mcom w_n)=\\
\sum
\begin{Sb}
k+l=n+1\\
\\
0\le i\le n-k
\end{Sb}
\s{(w_1\mplus w_i)\delta_k}
d_l(w_1\mcom w_i,\delta_k(w_{i+1}\mcom w_{i+k}),w_{i+k+1}\mcom w_n),
\end{multline}
for $d$, $\delta\in C(W)$.
Then $[d,\delta]=d\delta-\s{d\delta}\delta d$. One can express this in the
form 
\begin{equation}
[d,\delta]_n=
\sum
\begin{Sb}
k+l=n+1\\
\end{Sb}
[d_k,\delta_l],
\end{equation}
where $d_k$ represents the degree $k$ part of $d$, etc.
Since $d$ is odd, the fact that it is a codifferential can be expressed in the 
simple form $[d,d]=0$. The differential $D$ determining the cohomology
of the \ainf\ algebra is given by $D(\delta)=[\delta,d]$, and the fact
that $D^2=0$ follows from the (\zt-graded) Jacobi identity. The homology
of $C(W)$, or the corresponding homology of $C(V)$, is called the 
Hochschild cohomology of the \ainf\ algebra.

Now suppose that $V$ is equipped with a graded symmetric inner product,
that is a non-degenerate bilinear form  $h:V^2\ra\k$, satisfying
$h(v_1,v_2)=\s{v_1v_2}h(v_2,v_1)$. (Non-degeneracy means that the
induced map $h':V\ra V^*$, given by $h'(v_1)(v_2)=h(v_1,v_2)$ 
is an isomorphism.) Denote $h(v_1,v_2)=\ip{v_1}{v_2}$ for simplicity.
The inner product induces an isomorphism
between $C^n(V)=\hom(V^n,V)$ and $C^{n+1}(V,\k)=\hom(V^{n+1},\k)$,
given by $\ph\mapsto\tph$, where
\begin{equation}\label{cyceq}
\tph(v_1\mcom v_{n+1})=\ip{\ph(v_1\mcom v_n)}{v_{n+1}}.
\end{equation}
$\ph_n\in C^n(V)$ is said to be 
cyclic with respect to the inner product if
\begin{equation}\label{cycdef}
\ip{\ph_n(v_1\mcom v_n)}{v_{n+1}}=\s{n+v_1\mu}
\ip{v_1}{\ph_n( v_2\mcom v_{n+1})}.
\end{equation}
This is equivalent to the condition that $\tph_n$ is cyclic in the sense
that
\begin{equation}
\tph_n(v_1\mcom v_{n+1})=
\s{n+v_{n+1}(v_1\mplus v_n)}\tph_n(v_{n+1},v_1\mcom v_n).
\end{equation}
A cochain $\ph$ is said to be cyclic if $\ph_n$ is cyclic for all $n$.
Note that $\tph_n$ is cyclically symmetric for even $n$, but
cyclically antisymmetric for odd $n$. This feature contributed to
the difficulty in defining the homology cycle in \cite{ps2}. 
The inner product is said to be invariant with respect to the \ainf\
structure $m$ if $m$ is cyclic with respect to the inner product.
The \zt-graded bracket on $C(W)$ induces a \ztz-graded bracket 
$\brf$ on $C(V)$, and it was shown in \cite{pen1} that the bracket 
of two cyclic cochains is again cyclic. For an invariant inner product,
the cyclic elements form a subcomplex $CC(V)$ of $C(V)$, and the 
homology of this subcomplex is called the cyclic cohomology of the 
\ainf\ algebra.

Let us recast the notion of cyclic homology in terms of the complex
$C(W)$. We shall see that the definition becomes more natural. First
let us consider the inner product $k$ on $W$ induced by the inner product
$h$ on $V$.  This inner product is simply $h\circ\eta$, and it is easy
to check that
\begin{equation}
k(w_2,w_1)=\s{w_1w_2+1}k(w_1,w_2),
\end{equation}
so that $k$ is graded antisymmetric.

Now one can define cyclicity in $C(W)$ in one of three manners. We
could simply define $\psi$ to be cyclic if 
$\psi=\eta^{-1}\circ\ph\circ \eta$, for a cyclic cochain $\ph$ in 
$C(V)$. On the other hand, one could define $\psi$ to be cyclic
the corresponding element
$\tilde\psi\in C(V,\k)$ given by $\tph=\tilde\psi\circ\eta^{-1}$
is cyclic. It turns out that these notions are equivalent, and in
addition they are equivalent to the following definition of cyclicity.
We say that $\psi_n\in C^n(W)$ is cyclic with respect to the 
inner product $k$ if
\begin{equation}\label{cycpar}
\ip{\psi(w_1\mcom w_n)}{w_{n+1}}=\s{w_1\psi+1}
\ip{w_1}{\psi(w_2\mcom w_{n+1})}
\end{equation}
This notion is also equivalent to $\tilde\psi$ being cyclic in the
sense that
\begin{equation}
\tilde\psi(w_1\mcom w_{n+1})=\s{w_{n+1}(w_1\mplus w_n)}
\tilde\psi(w_{n+1},w_1\mcom w_n)
\end{equation}
Note that $\tilde\psi$ is cyclically symmetric for all $n$.
Next, we note that if $\delta$ and $\psi$ are cyclic cochains
in $C(W)$, then their bracket is cyclic, and in fact, we have the
following theorem.
\begin{thm} Suppose that $W$ is a \zt-graded vector space with
a graded antisymmetric inner product. 

a) If $\delta$, $\psi\in C(W)$
are cyclic with respect to the inner product, then their bracket
is also cyclic. Moreover, the following formula holds.
\begin{multline}
\widetilde{[\delta,\psi]}(v_1\mcom v_{n+1})=\\
\sum
\begin{Sb}
k+l=n+1
\\
0\le i\le n\\
\end{Sb}
\s{(w_1\mplus w_i)(w_{i+1}\mplus w_{n+1})}
\tilde\delta_l(\psi_k(w_{i+1}\mcom w_{i+k}),w_{i+k+1}\mcom w_i).
\end{multline}
Thus there is a bracket defined on the complex $CC(W,\k)$ of cyclic
elements in $C(W,\k)$, by 
$[\tilde\delta,\tilde\psi]=\widetilde{[\delta,\psi]}$.
b) If $d$ is a \zt-graded codifferential on $T(W)$, then there is
a differential $D$ in $CC(W,\k)$, given by 
\begin{multline}
D(\tilde\psi)(w_1\mcom w_n)=\\
\sum
\begin{Sb}
k+l=n+1
\\
0\le i\le n\\
\end{Sb}
\s{(w_1\mplus w_i)(w_{i+1}\mplus w_{n+1})}
\tilde\psi_l(d_k(w_{i+1}\mcom w_{i+k}),w_{i+k+1}\mcom w_i).
\end{multline}

c) If the inner product is invariant, then 
$D(\tilde\psi)=[\tilde\psi,\td]$. Thus $CC(W,\k)$ inherits the 
structure of a differential graded Lie algebra.
\end{thm}

This result is really just a restatement of the corresponding result
for cyclic cohomology formulated in terms of $C(V)$ (see cite{pen1}). 
The form in which we will
need to use this result is as follows. Suppose that $e_1\mcom e_n$ is
a basis of $W$. (We assume always that $V$ is finite dimensional.)
The structure constants for the algebra are given by 
$d_k(e_{j_1}\mcom e_{j_k})=d^a_{j_1\mcom j_k}e_a$. (We use the summation
convention for repeated upper and lower indices.)
Define the lower structure constants 
$d_{j_1\mcom j_{k+1}}$ of the \ainf\ algebra
by 
\begin{equation}\label{structrel}
d_{j_1\mcom j_{k+1}}=\td_k(e_{j_1}\mcom e_{j_{k+1}}).
\end{equation}
If the inner product is denoted by 
$k_{\alpha\beta}=\ip{e_\alpha}{e_\beta}$, then $k^{\alpha\beta}$ denotes
the inverse matrix to $k$. Then the lower structure constants are
given in terms of the usual ones by
$d^a_{j_1\mcom j_k}e_a=h^{a,j_{k+1}}d_{j_1\mcom j_{k+1}}$.
Then we can formulate the definition of an \ainf\ algebra with an
invariant inner product in terms of the structure constants by
\begin{equation}\label{ainfrel}
\sum
\begin{Sb}
k+l=n+1
\\
0\le i\le n\\
\end{Sb}
\s{(e_{j_1}\mplus e_{j_i})(e_{j_{i+1}}\mplus e_{i_{n+1}})}
d^a_{j_{i+1}\mcom j_{i+l}}d_{a,j_{i+l+1}\mcom j_{i}}=0.
\end{equation}
This equation will play an important role in the construction of
the cycle in the homology of the complex of metric ribbon graphs.
Let us also summarize the additional information that we shall 
need. Although we did not state this before, we shall need to
assume that the inner product is an even map. It is also antisymmetric,
so that the tensor $K=k^{ab}e_a\tns e_b$ 
is an even, antisymmetric tensor. 
The tensor $C_k=d_{j_1\mcom j_{k+1}}e^{j_1}\mtns e^{j_k}$ is an odd,
cyclically symmetric tensor.
\section{\linf\ algebras}\label{linfsect}
In this section we again assume that $V$ is a finite dimensional
\zt-graded vector space over a field \k of characteristic zero.
The exterior coalgebra $\bigwedge V$ is defined by 
$\bigwedge V=\bigoplus_{n=1}^\infty \bigwedge^n V$, where $\bigwedge ^n V$
is the $n$-th exterior power of $V$. An \linf\ structure on $V$
is given by an element (cochain) $l\in the complex
C(V)=\hom(\bigwedge V,V)$, so that $l$ is a collection of maps
$l_k:\bigwedge V\ra V$. We require that the relations
\begin{equation}
\sum
\begin{Sb}
k+l=n+1
\\
\sigma\in\sh{l}{n-l}\\
\end{Sb}
\s{\sigma}\epsilon(\sigma)\s{(k-l)l}
l_l(l_k(v_{\sigma(1)}\mcom v_{\sigma(k)})
,v_{\sigma(k+1)}\mcom v_{\sigma(n)}=0.
\end{equation}
where $\sh{l}{n-l}$ is the set of unshuffles of type $(l,n-l)$,
in other words, permutations $\sigma$ satisfying 
$\sigma(i)\le \sigma(i+1)$ unless $i=l$, $\s{\sigma}$ is the
sign of the permutation, and $\epsilon(\sigma)$ is a sign  which
depends on both $\sigma$ and $v_1\mcom v_n$. For a more detailed 
explanation of these signs, see \cite{pen1}. In addition, the maps
$l_k$ are of parity $k$. In complete parallel to the \ainf\ picture,
we have $l_1$ is an odd differential on $V$, acting as a derivation
on the multiplication $l_2$, which satisfies the Jacobi identity up
to a homotopy determined by $l_3$. \linf\ algebras are also called
{\bf strongly homotopy Lie algebras}. Lie algebras and differential
graded Lie algebras are special cases of \linf\ algebras.

Let $\Pi W=V$ be the parity reversion of $V$. There is a natural isomorphism
$\eta:\bigodot W\ra \bigwedge V$, 
where $\bigodot W$ is the symmetric coalgebra of $W$, which
is defined by 
\begin{equation}
\eta(w_1\modot w_n)=
\s{(n-1)w_1\mplus w_{n-1}}
\pi w_1\mwedge w_{n}.
\end{equation}
If we denote $C(W)=\hom(\bigodot W,W)$, then the isomorphism $\eta$ 
induces an isomorphism between $C(V)$ and $C(W)$ given by 
$\ph\mapsto \eta^{-1}\circ\ph\circ\eta$. It is well known 
that $l$ determines an \ainf\ structure precisely when its image
$d$ under this isomorphism determines a \zt-graded odd codifferential
on $\bigodot W$. We shall regard $W$ equipped with $d$ as an \linf\
algebra. The statement that $d$ determines a  codifferential
is equivalent to
\begin{equation}
\sum
\begin{Sb}
k+l=n+1
\\
\sigma\in\sh{l}{n-l}\\
\end{Sb}
\s{\sigma}\epsilon(\sigma)
d_l(d_k(v_{\sigma(1)}\mcom v_{\sigma(k)}),
v_{\sigma(k+1)}\mcom v_{\sigma(n)}=0.
\end{equation}
Elements of $C(W)$ extend uniquely to coderivations on $\bigodot W$,
so $C(W)$ has a \zt-graded Lie bracket, and the codifferential  $d$
determines a differential $D$ on $C(W)$ by $D(\psi)=[\psi,d]$. This
differential determines the homology of the \linf\ algebra. 
Cocycles
of a Lie algebra give rise to infinitesimal deformations of the Lie
algebra into an \linf\ algebra.

Now suppose $V$ has a \zt-graded symmetric inner product $h$, 
so that $W$ is equipped with a \zt-graded antisymmetric inner product
$k$. An element $\ph\in\hom(\bigwedge^n V,V)$ is cyclic if
it satisfies equation (\ref{cycdef}). The map $\ph$ is cyclic if and
only if the map $\tilde\ph$, given by equation (\ref{cyceq}) is 
antisymmetric, in other words, $\tilde\ph\in C^{n+1}(V,\k)=\hom(
\bigwedge^{n+1} V,\k)$. Denote the set of cyclic elements by $CC(V)$.
The Lie bracket on $C(W)$ induces a bracket on $C(V)$, and the bracket
of two cyclic elements in $C(V)$ is again cyclic, so $CC(V)$ is a Lie
subalgebra of $C(V)$. If the inner product is invariant, then $CC(V)$
is a subcomplex of $C(V)$. The cyclic cohomology of the \linf\ is the
homology of this subcomplex.

One also can define cyclic cohomology using the complex $C(W)$, and 
the antisymmetric inner product $k$. The  definition of cyclicity
of an element $\psi\in C(W)$ is again given by equation (\ref{cycpar}).
The notion of cyclicity is equivalent to the map $\tilde\psi$, given
by 
$\tilde\psi(w_1\mcom w_{n+1})=\ip{\psi(w_1\mcom w_n)}{w_{n+1}}$,
being graded symmetric.
The bracket of cyclic elements in $C(W)$ is again cyclic, and so we
obtain the following theorem.
\begin{thm} Suppose that $W$ is a \zt-graded vector space with
a graded antisymmetric inner product. Let $C(W)=\hom(\bigodot W,W)$
be equipped with the bracket induced by the identification of elements
in $C(W)$ with coderivations of $\bigodot W$.

a) If $\ph$, $\psi\in C(W)$
are cyclic with respect to the inner product, then their bracket
is also cyclic. Moreover, the following formula holds.
\begin{multline}
\widetilde{[\ph,\psi]}(v_1\mcom v_{n+1})=\\
\sum
\begin{Sb}
k+l=n+1
\\
\sigma\in\sh{l}{n+1-l}\\
\end{Sb}
\s{\sigma}\epsilon(\sigma)
\tilde\ph_l(\psi_k(w_{\sigma(1)}\mcom w_{\sigma(k)}),
w_{\sigma(k+1)}\mcom w_{\sigma(n+1)}).
\end{multline}
Thus there is a bracket defined on the complex $CC(W,\k)$ of cyclic
elements in $C(W,\k)$, by
$[\tilde\ph,\tilde\psi]=\widetilde{[\ph,\psi]}$.

b) If $d$ is a \zt-graded codifferential on $T(W)$, then there is
a differential $D$ in $CC(W,\k)$, given by
\begin{multline}
D(\tilde\psi)(w_1\mcom w_{n+1})=\\
\sum
\begin{Sb}
k+l=n+1
\\
\sigma\in\sh{l}{n+1-l}\\
\end{Sb}
\s{\sigma}\epsilon(\sigma)
\tilde\ph_l(d_k(w_{\sigma(1)}\mcom w_{\sigma(k)}),w_{\sigma(k+1)}\mcom
w_{\sigma(n+1)})
\end{multline}

c) If the inner product is invariant, then
$D(\tilde\psi)=[\tilde\psi,\td]$. Thus $CC(W,\k)$ inherits the
structure of a differential graded Lie algebra.
\end{thm}
The form in which we need these results is as follows. Let $e_1\mcom e_m$
be a basis of $W$, and define the lower structure constants of the \linf\
algebra $d_{j_1\mcom j_k}$ by equation (\ref{structrel}). Then we can
reformulate the definition of an \linf\ algebra with an antisymmetric
invariant inner product in terms of the structure constants by
\begin{equation}\label{linfstruc}
\sum
\begin{Sb}
k+l=n+1
\\
\sigma\in\sh{l}{n-l}\\
\end{Sb}
\epsilon(\sigma)
d^a_{j_{\sigma(1)}\mcom j_{\sigma(l)}}
d_{a,j_{\sigma(l+1)}\mcom j_{\sigma(n)}}=0.
\end{equation}
\section{Definition of the Graph Complexes}
A graph is a 1 dimensional CW complex, in other words, it consists of
vertices and edges. We will consider here only graphs where each 
vertex is at least trivalent, meaning that the number of edges incident
to the vertex is at least 3. The graph is called a ribbon graph if
in addition there is a fixed cyclic order of the edges at each vertex.
A metric on the graph is an assignment of a positive number to each
edge of the graph. The set $\sg$ of all metrics on a graph
gives a cell, identifiable with $R_+^k$ where $k=e(\Gamma)$ is the
number of edges in the graph. Intuitively, when the length of an
edge of a graph tends to zero, the graph degenerates to a new graph
$\Gamma'$, which determines a cell on the boundary of the cell 
$\sg$. We shall give a more precise definition of the boundary 
operator shortly. Let us point out that neither 
the space of metric ribbon graphs, nor that of metric ordinary graphs
is a cell complex, because the closure of a cell is not compact. (Boundary
cells exist only on the sides corresponding to setting a length equal to
zero, not on the infinite sides.)

If we consider the space of ordinary graphs, then the Euler Characteristic
is preserved when contracting edges, so there are different graph complexes
corresponding to the different Euler characteristics.  For the space of
metric ribbon graphs, each metric graph corresponds to a Riemann surface
with a fixed genus $g$ and number of marked points $n$, depending only
on the graph (see \cite{har,pen}). 
The boundary operator also preserves the genus and number
of marked points. Thus for ribbon graphs, we
have different graph complexes for each
genus and number of marked points. 

To define a boundary operator, we shall need a notion of the orientation
of the cells in the graph complex. 
Let us first proceed with the case
of ribbon graphs. An orientation of the complex of metric ribbon graphs
is determined by a ordering of the edges in the graph, and a ordering
of the holes or circuits in the graph (see \cite{kon,ps2}). Using this
notion of orientation, a cycle in the complex of metric ribbon graphs
associated to an \ainf\ algebra with a symmetric invariant inner product
was constructed in \cite{ps2}. This construction was dependent on the
non trivial fact that the complex of metric ribbon graphs is orientable.
Here we propose an equivalent definition of the orientation of a cell
corresponding to a graph $\Gamma$. The orientation of the cell is given
by choosing a ordering of the vertices and assigning an arrow to each
edge (in other words, choosing an orientation of the edge). The fact that
this orientation is equivalent to the previous one is not obvious, so 
here we also are concealing some non trivial relations. For the
case of metric ordinary graphs, in  \cite{kon3} an orientation was
defined to be an ordering of the edges in the graph, as well as an
orientation of the first homology group. It is not difficult in this
case to see that this notion of orientation is equivalent to the definition
of an orientation as an ordering of the vertices and an assignment of
an orientation to each edge of the graph. Thus our definition of 
orientation is the same for ribbon graph and ordinary graph complexes.

Next, we would like to define the boundary operator.
Let us make the abuse of notation and identify the cell $\sg$ 
with the graph $\Gamma$. Two oriented graphs are considered to be
the same graph if there is an isomorphism between the graphs which
induces the same orientation. If we contract out an edge in an
oriented graph, then there is a natural orientation induced in the
contracted graph from an orientation in the original graph. To see
this, suppose that the graph $\go$ has edges labeled from 1 to $n$,
and vertices labeled from 1 to $m$. When an edge is contracted out,
a vertex also is contracted. An orientation only determines the
labels on the edges and vertices up to an even permutation, so it is
always possible to assume that the labeling has been chosen so that
the edge to be contracted has a vertex labeled as $m$, and that the
arrow points towards this vertex. Then, when contracting, the vertex
with label $m$ disappears. However, the process is reversible, since
when inserting an edge, one has a choice of which of the two vertices
that occur will have the new label, requiring that the arrow
in the inserted edge point to the new vertex fixes the orientation
uniquely. Notice that the same construction applies to ordinary and
to ribbon graphs. The only special feature for ribbon graphs is that
when contracting edges, in combining the vertices, there is an induced
cyclic order on the vertex, and similarly, when inserting an edge,
there is a natural way to induce a cyclic order at the new vertices
from the cyclic order on the old vertex. 

The graph complex is generated by the classes of oriented
graphs, with relations given by identifying an oriented 
graph with the negative of the same graph with opposite orientation.
In our construction, we are taking coefficients in the field \k,
which has characteristic zero, so that if an oriented
 graph is equivalent to
its opposite orientation, then it is zero in the graph complex. 
But it is also interesting to consider the case of coefficients in \zt,
in which case orientation plays no role, so that one obtains a different
homology theory.
A graph with a loop is equivalent to its opposite orientation. To see this,
note that an equivalence of oriented graphs is given by an equivalence
of graphs, which assigns a vertex to each vertex, and thus induces an
ordering of the vertices. Normally, one can assign an arrow to an arrow,
by requiring that the induced arrow on an edge point from the induced 
starting vertex to the induced ending vertex, but this procedure is 
ambiguous when the starting and ending vertex is the same. Thus either
orientation can occur as the image of the oriented graph.

Let us introduce a boundary operator in the following manner. If 
$\go$ and $\gop$ are oriented graphs, then we need to assign an incidence
number $[\gop,\go]$, which counts the number of times that $\gop$ occurs
in the boundary of $\go$. 
Let $G$ be the set of all distinct equivalence classes in the complex.
Graphs with opposite orientation determine the same equivalence class
in this context.
In the case of ribbon graphs, we can choose $G$ to be the set of all 
equivalences of some fixed genus $g$ and number of marked points $n$,
while for the ordinary graph complex, we can choose $G$ to be the set
of all equivalence classes of the same Euler characteristic.
A chain in the complex can be represented
as a sum of the form 
$\sum_{\go\in G} a_{\go}\go$, 
Then the boundary operator is given
by 
\begin{equation}
\del(\go)=
\sum_{\go\in G}
[\gop,\go]\gop
\end{equation}
To determine $[\gop,\go]$, we consider the result of contracting an edge
in $\go$. An edge is contractable if it is attached to two distinct vertices.
When an edge is contracted, the graph $\gop$ or its opposite may occur, and
if so, one counts either a plus one or a minus one, depending on which occurs.
The sum of the numbers resulting from contracting the
various edges in \go\ is the desired incidence number.

\begin{figure}[thb]
\epsfxsize=2.4 in
\centerline{\epsfbox{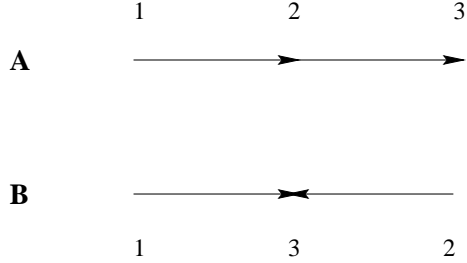}}
\caption{Contracting adjacent edges}\label{boundary}
\end{figure}

\begin{figure}[thb]
\epsfxsize=2.4 in
\centerline{\epsfbox{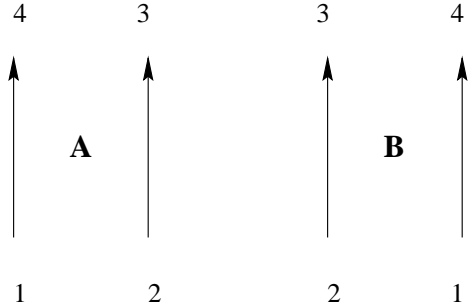}}
\caption{Contracting non adjacent edges}\label{boundfig}
\end{figure}

Figures \ref{boundary} and
 \ref{boundfig} illustrate why the square of the boundary operator
is zero. When the boundary operator is applied twice, that is the same
as contracting two edges in the graph. These edges may be contracted in
two orders, and the resulting graphs will be the same. Thus we need only
show that the incidence numbers associated to these two different orders
are opposite, to show that the square of the boundary operator is zero.
In figure \ref{boundary}
the orientations in A and B are the same. Recall that to contract
an edge, one should place the highest number of a vertex at one of the
edges, and the arrow on the edge should face that vertex. Thus in A, we
can contract the edge from vertex 1 to 3, and then the arrows will be in
accord so that the remaining edge can be contracted. However, in B, we
can contract the edge from 3 to 2, but then the remaining arrow is 
pointing in the wrong direction. This shows that the consequence of 
contracting the two edges in the opposite order is that the incidence
numbers will cancel.
\vskip .10 in
\vskip .10 in
In figure \ref{boundary}, the two edges were assumed to be adjacent. If
the edges are not adjacent, one can arrange the labels as in figure
\ref{boundfig} below. The orientations in A and B below are the same. In A,
one can contract first the edge with vertices labeled 1 and 4, and then
the edge with vertices labeled 3 and 2. The same procedure is applied to
B. The resulting graphs have opposite orientation, since the position of
the vertices 1 and 2, which remain after the contraction, is reversed.
\vskip .10 in
\vskip .10 in
\section{The cycle given by an \ainf\ algebra}
Let $W$ be finite dimensional vector space with basis $e_i$, equipped
with an \ainf\ algebra structure given by a codifferential $d$ on $T(W)$.
Suppose that $k$ is a graded antisymmetric inner product on $W$ which is
invariant with respect to the \ainf\ structure, and that the tensors
$K=k^{ab}e_a\tns e_b$ and $C_n=d_{j_1\mcom j_n}e^{j_1}\mtns e^{j_n}$ are
defined as in section \ref{ainfsect}. Recall that $K$ is antisymmetric and
even, and that $C_n$ is odd and cyclically symmetric. 
We want to construct a cochain in the graph complex, which is a function
$Z$ on the set of oriented graphs such that $Z(\go)$ has opposite sign
for graphs of opposite orientation. To define $Z$, let $\go$ be an
oriented graph. and to each edge in the graph associate two indices, 
one for each incident vertex. 
Intuitively, if the edge has labels $i$ and $j$, then
associate the tensor $h^{ij}$, and if the labels on a vertex with
$n$ incident edges are
$j_1\mcom j_n$ in cyclic order, then associate the tensor
$d_{j_1\mcom j_n}$ to the vertex. Multiply the tensors for the vertices
in the same order as the ordering of the vertices, and then multiply
the tensors corresponding to the edges in any order. Then one obtains
a tensor in which every upper indice matches precisely one lower index,
so summing over repeated indices yields a number Z(\go). 
Figure \ref{contract}
below illustrates the idea. The vertices are labeled 1 and 2. The 
cyclic order at each vertex is assumed to be counterclockwise.
The portion of the product corresponding to the figure is
$d_{glmih}d_{nkj}k^{bh}k^{ic}k^{mn}k^{lf}k^{ag}k^{jd}k^{ek}$.
\vskip .10 in
\begin{figure}[htb]
\epsfxsize=2.4 in
\centerline{\epsfbox{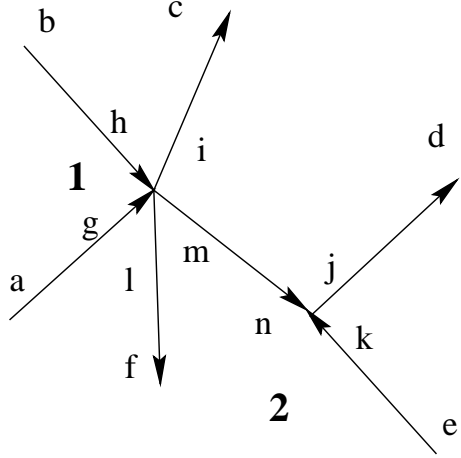}}
\caption{Constructing the cycle}\label{contract}
\end{figure}
\vskip .10 in
More precisely, one should consider the product of tensors $C_n$, in the
order dictated by the order of the vertices, and the tensors $K$, one for
each edge. This will yield an element of $(W^*)^{2e}\tns W^{2e}$. Then
one uses the graph to determine a graded contraction of this tensor,
obtaining an element of \k. In the informal construction, one must take
into account a sign that arises from the graded contraction as well. The
sign corresponding to the above portion arises as a consequence of the
exchange rule in computing the contraction. For example,
if one contracts the tensor
$e^a\tns e^b\tns e_a\tns e_b$ according to the order indicated, then
one obtains the sign $\s{e^be_a}$, but the contraction of the tensor
$e^a\tns e^b\tns e_b\tns e_a$ yields the sign 1.

We need to show that the definition of the number $Z(\go)$ does not depend
on any of the choices that we made in writing down the tensor. First of
all, since the tensors $C_n$ are odd, permuting their order will introduce
a sign corresponding exactly to the sign of the permutation, so this is
precisely what we need to have the sign depend properly on the orientation.
Secondly, $C_n$ is cyclically (graded) symmetric, so that the starting point
we used to write the tensor down makes no difference. Thirdly, the order
in which the tensors $K$ are written down makes no difference, since they
are even tensors. Finally, the order in which the indices are written down
must conform to the order in which the arrow appears, and the antisymmetry
guarantees that the reversal of the arrow will reverse the sign in the 
contraction. Thus we have shown that 
\begin{equation}
Z=\sum_{\go\in G}Z(\go)\go
\end{equation}
is a well defined chain in the complex.

Next, we wish to show that $Z$ is a cycle. It is sufficient to consider
a graph $\gop$, and consider a single vertex in $\gop$, and the graphs
which are obtained by inserting an edge in the vertex. We show that the
sum of the contributions of each of these graphs to the boundary of $Z$
is zero. Let us consider how an  edge might be inserted in a vertex with
$n+1$ incident edges to obtain two new vertices with $k+1$ and $l+1$ 
incident edges, where $k+l=n+1$,
as is illustrated in the figure below. 
\vskip .10 in
\begin{figure}[htb]
\epsfxsize=4 in
\centerline{\epsfbox{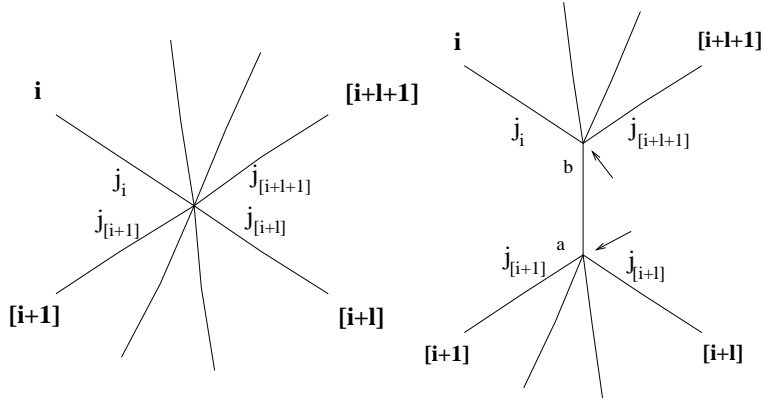}}
\caption{Splitting a vertex}\label{versplit}
\end{figure}
\vskip .10 in
When inserting a new edge, one
specifies two incident edges, $j_i$ and $j_{i+l}$ in the diagram. If the
cyclic order at the vertex is counterclockwise, then the cyclic order
at the two new vertices is counterclockwise. Denote the expanded graph
by \go.
Let us suppose that the vertex with new incident edge labeled $b$ is the
one which is given the new highest label. Then the arrow runs from 
$a$ to $b$, and the contribution to the function $Z(\go)$ from these two
vertices and the new edge is 
$d_{
j_{i+1}\mcom j_{i+l},a}
d_{b,j_{i+l+1}\mcom j_{i}}h^{ab}$. 
Since $d^b_{j_{i+1}\mcom j_{i+1}}=d_{j_{i+1}\mcom j_{i+l},a}h^{ab}$,
and the sign corresponding to the graded contraction is 
$\s{(e_{j_1}\mplus e_{j_i})(e_{j_{i+1}}\mplus e_j{n+1})}$. we see that
this number is just the $i$-th term in the expression given in 
equation \ref{ainfrel}. Since the sum of all the terms is zero.
we conclude that the contribution of all the graphs which arise from
the insertion of an edge in \gop\ is zero. Thus $Z$ is a cycle.
In the proof, we are tacitly assuming that the $d_1$ term in the
\ainf\ algebra vanishes, but the proof can be extended to show that
the chain $Z$ is still a cycle when the $d_1$ term does not vanish,
see \cite{pen2}.

In \cite{ps2}, we showed that if $A$ was an associative algebra
with product $m$, and 
$\ph$ was a $k$-cocycle, then one obtains a cycle in the homology
of the graph complex by considering the infinitesimal deformation
$m_t=m+t\ph$ of the associative algebra into an \ainf\ algebra.
The cycle so constructed depends only on graphs which 
have trivalent vertices, with the exception of one vertex which is 
$k$-valent. One can prove this directly as well, using the corresponding
formula for the cyclic coboundary operator on the parity reversion.
Albert Schwarz conjectured that the cycle obtained in this manner was
independent of the algebra in the sense that one always obtains a 
multiple of the same cycle. As it turns out, this conjecture follows from
a theorem of Penner in \cite{pnr2}, because there is only one cycle of
this type.
\section{The cycle given by an \linf\ algebra}
The cycle in the graph complex of metric ordinary graphs is constructed in 
the same manner as the cycle in the complex of ribbon graphs, so I shall
just describe the proof that it is a cycle, which is a little different.
For ordinary graphs, there is no cyclic order at each vertex, so the
number of ways in which an edge can be inserted at a vertex is much 
larger. For example, if a vertex has four incident edges, labeled 
$a$, $b$, $c$ and $d$, then if the graph is a ribbon graph, and the
cyclic order is $a,b,c,d$, one can insert an edge so that $a$ and $b$
are paired with one of the vertices, and $c$ and $d$ with the other.
Alternatively, one can insert an edge so that $b$ and $c$ pair with
one of the vertices, and $d$ and $a$ with the other. No other possibility
occurs. However, for the ordinary graph complex, we cannot distinguish
any order, so there is in addition, the possibility that $a$ and $c$ pair
off with a new vertex, and $b$ and $d$ with the other. The combinations
that can occur when decomposing a vertex with $n+1$ incident edges into
two vertices with $l+1$ and $k+1$ edges respectively are given by the
shuffles of type $(l,k)$. 
Let \gop\ be a graph and consider a vertex with $n+1$ incident edges
labeled $j_1\mcom j_{n+1}$.
If $\sigma$ is such a shuffle, then let \go\ be the graph 
which results from inserting
an edge so that the edges $\sigma(1)\mcom \sigma(l)$ are associated to
one of the new vertices, and the other edges with the other vertex.
Then the contribution to $Z(\go)$ from the new vertices and the edge
is $d_{j_{\sigma(1)}\mcom j_{\sigma(l)},a}
d_{b,j_{\sigma(l+1)}\mcom j_{\sigma(n+1)}}$. The sign corresponding to
the graded contraction turns out to be just $\epsilon(\sigma)$. Then
it is clear from equation (\ref{linfstruc}) that the sum of the terms
corresponding to the various ways to insert the new edge cancel in their
contribution to the boundary of $Z$.

In \cite{pen1}, it was shown that if $L$ is a Lie algebra with bracket
$l$ and invariant inner product, then a cyclic cocycle $\ph$ determines an
infinitesimal deformation of the Lie algebra into an \linf\ algebra
with invariant inner product. Suppose that the cocycle is of exterior
degree $n$, in other words, $\ph\in C^n(V)$. Then we can associate to
this cocycle a cycle in the homology of the ordinary graph complex 
which depends only on the graphs with all vertices trivalent, except
for one vertex with $n+1$ edges. I do not know if the cycle depends
on the algebra.

\section{Conclusion}
This paper has shown how to construct cycles in the homology of graph
complexes from \ainf\ and \linf\ algebras with invariant inner products.
A corollary of this construction is the existence of cycles in the
graph complexes associated to cyclic cocycles of associative and Lie
algebras with invariant inner products. These cycles exist because
of the connection between the cyclic cohomology of these algebras
and infinitesimal deformations of the algebras into infinity algebras.
One can also consider higher order deformations, leading to the 
construction of some other cycles in the homology of the graph complexes.

There is a general theory of homotopy algebras, which gives \ainf\
and \linf\ algebras as special cases. This theory was introduced in
\cite{getz}
In addition, there is a unified
theory, called operad cohomology, which was introduced by T. Fox and
M. Markl in \cite{fm}. The graph complexes are closely related to
certain operad structures. There must be some relation between the 
results here and the operads which are used to define these homotopy
algebras. Probably a simple proof of these results can be given,
using operad theory.

One can generalize the graph complexes to allow vertices with only
two edges. In this more general complex, it is immediate that one
can define cycles using \linf\ and \ainf\ terms which do not have
trivial $d_1$ term. However, it was shown in \cite{pen2}, and is
not difficult to see, that the generalized graph complex is a
direct sum of the ordinary graph complex and another subcomplex.
As a consequence, it follows that \linf\ and \ainf\ algebras
with non trivial $d_1$ term give rise to cycles in the homology
of the graph complex. One can also take the homology of an 
infinity algebra with respect to the $d_1$ term, and this homology
inherits the structure of an infinity algebra with a trivial $d_1$
term. The cycle given by this homology algebra probably coincides
with the cycle given by the original algebra, but I don't know of
a proof of this result.
\bibliographystyle{amsplain}
\ifx\undefined\bysame
\newcommand{\bysame}{\leavevmode\hbox to3em{\hrulefill}\,}
\fi

\end{document}